\documentclass[12pt]{amsart}
\usepackage{amsmath, amssymb, latexsym, amsthm, mathrsfs, geometry}
\usepackage[all]{xy}
\usepackage{color}

\newcommand{\nc}[1]{\newcommand{#1}}
\nc{\ed}{
\subsection*{Acknowledgments}
We thank Angelo Bella for useful discussions and for bringing Aurichi's results to our attention.
We thank Leandro Aurichi for making his preprint \cite{Aurichi} available to us. Proposition
\ref{Rcoc} is based on one by Aurichi.
A part of this research was carried out during a visit of the third named author in the University of Messina.
This author thanks his hosts, and in particular the other named authors and Andrei Catalioto,
for their kind hospitality.

\end{document}}

\nc{\D}{\mathrm{D}}\nc{\B}{\mathrm{B}}
\nc{\Do}{\D_\mathrm{o}}
\nc{\my}[1]{\textcolor{green}{\sf #1}}
\nc{\Pa}[8]{\bibitem{#1} {#2}, \emph{#3}, {#4} \textbf{#5} ({#6}), {#7}--{#8}.}
\nc{\tPa}[5]{\bibitem{#1} {#2}, \emph{#3}, {#4}, to appear.}
\nc{\sPa}[4]{\bibitem{#1} {#2}, \emph{#3}, {#4}, submitted.}
\nc{\Bc}[9]{\bibitem{#1} {#2}, \emph{#3}, in: \textbf{#4} (#5), #6 #7, #8--#9.}
\nc{\fD}{\mathfrak{D}}\nc{\fX}{\mathfrak{X}}\nc{\Dfin}{\mathfrak{D}_\mathrm{fin}}
\nc{\ft}{\mathfrak{t}}\nc{\fb}{\mathfrak{b}}\nc{\fc}{\mathfrak{c}}\nc{\fd}{\mathfrak{d}}
\nc{\fg}{\mathfrak{g}}\nc{\fr}{\mathfrak{r}}\nc{\fu}{\mathfrak{u}}\nc{\fh}{\mathfrak{h}}
\nc{\fp}{\mathfrak{p}}\nc{\fj}{\mathfrak{j}}\nc{\fs}{\mathfrak{s}}\nc{\od}{\mathfrak{od}}
\nc{\Onbd}{\Op_{\mathrm{nbd}}}\nc{\Omnb}{\Omega_{\mathrm{nbd}}}
\nc{\pr}{\op{pr}}\nc{\plim}{p\txt{-}\lim}\nc{\osc}{\op{osc}}\nc{\supp}{\op{supp}}\nc{\PR}{\op{PR}}
\nc{\Bgp}{{\Z^\N}}\nc{\Cgp}{{{\Z_2}^\N}}\nc{\Cantor}{{\{0,1\}^\N}}
\nc{\Fr}{\mathit{F\!r}}
\nc{\Cite}[1]{\textbf{[#1]}}
\nc{\Next}[1]{{#1^+}}\nc{\intvl}[2]{{[#1(#2),\allowbreak #1(#2\!+\!1))}}
\nc{\gp}{\mbox{-\textit{\tiny gp}}}\nc{\grbl}{{\mbox{\textit{\tiny gp}}}}
\nc{\BOfat}{\B_{\Omega_{\mathrm{fat}}}}
\nc{\blocks}[2]{\op{cl}_{#2}(#1)}\nc{\blocksplus}[2]{\op{cl}^+_{#2}(#1)}
\nc{\arx}[1]{\texttt{http://arxiv.org/math/#1}}
\nc{\bq}{\begin{quote}}\nc{\eq}{\end{quote}}
\nc{\cl}[1]{\overline{#1}}
\nc{\CH}{the Continuum Hypothesis}\nc{\MA}{Martin's Axiom}
\nc{\inv}{^{-1}}
\nc{\bof}{\op{\fb}}\nc{\bofF}{\bof(\cF)}
\nc{\N}{\mathbb{N}}\nc{\bbP}{\mathbb{P}}\nc{\bbT}{\mathbb{T}}\nc{\Q}{\mathbb{Q}}\nc{\R}{\mathbb{R}}\nc{\Z}{\mathbb{Z}}
\nc{\NN}{{\N^{\N}}}\nc{\ZN}{{\Z^{\N}}}\nc{\NNup}{{\N^{\uparrow\N}}}\nc{\PN}{{P(\N)}}
\nc{\roth}{{[\N]^{\oo}}}\nc{\Fin}{{[\N]^{<\oo}}}\nc{\ici}{{[\N]^{(\oo,\oo)}}}
\nc{\compactN}{\N\cup\{\oo\}}\nc{\Inc}{{(\compactN)^{\uparrow\N}}}
\nc{\powInc}[1]{{\big(\Inc\big)^{#1}}}\nc{\powFin}[1]{{\big(\Fin\big)^{#1}}}\nc{\powPN}[1]{{\big(\PN\big)^{#1}}}
\nc{\NcompactN}{{(\compactN)^\N}}
\nc{\set}[2]{\{#1\,:\,#2\}}\nc{\seq}[1]{\{#1\}_{n\in\N}}\nc{\setseq}[1]{\set{#1}{n\in\N}}\nc{\sseq}[1]{\set{#1}{n\in\N}}
\nc{\op}{\operatorname}\nc{\im}{\op{im}}\nc{\maxfin}{\op{maxfin}}\nc{\ran}{\op{range}}
\nc{\scrA}{\mathscr{A}}\nc{\scrB}{\mathscr{B}}\nc{\scrC}{\mathscr{C}}\nc{\scrD}{\mathscr{D}}
\nc{\cA}{\mathcal{A}}\nc{\cK}{\mathcal{K}}\nc{\cD}{\mathcal{D}}\nc{\cF}{\mathcal{F}}\nc{\cS}{\mathcal{S}}
\nc{\cG}{\mathcal{G}}\nc{\cY}{\mathcal{Y}}\nc{\cL}{\mathcal{L}}\nc{\cM}{\mathcal{M}}\nc{\cN}{\mathcal{N}}
\nc{\cO}{\mathcal{O}}\nc{\cB}{\mathcal{B}}\nc{\cI}{\mathcal{I}}\nc{\cJ}{\mathcal{J}}
\nc{\bP}{\mathbf{P}}\nc{\bB}{\mathbf{B}}\nc{\bD}{\mathbf{D}}
\nc{\BG}{\B_\Ga}\nc{\BL}{\B_\Lambda}\nc{\BT}{\B_\Tau}\nc{\BTstar}{\B_{\Tau^*}}\nc{\BO}{\B_\Omega}
\nc{\BLgp}{\B_{\Lambda^{\grbl}}}
\nc{\CG}{\rmC_\Ga}\nc{\CL}{\rmC_\Lambda}\nc{\CT}{\rmC_\Tau}\nc{\CTstar}{\rmC_{\Tau^*}}\nc{\CO}{\rmC_\Omega}
\nc{\COgp}{C_{\Omega^{\grbl}}}\nc{\CLgp}{C_{\Lambda^{\grbl}}}\nc{\BOgp}{\B_{\Omega}^{\grbl}}
\nc{\sfC}{\mathsf{C}}\nc{\sfD}{\mathsf{D}}
\nc{\Tau}{\mathrm{T}}\nc{\Op}{\mathrm{O}}\nc{\Ga}{\Gamma}\nc{\Om}{\Omega}
\nc{\cP}{\mathcal{P}}\nc{\cU}{\mathcal{U}}\nc{\cV}{\mathcal{V}}\nc{\cW}{\mathcal{W}}
\long\def\forget#1\forgotten{}
\nc{\oo}{\infty}\nc{\w}{\omega}\nc{\x}{\times}\nc{\Iff}{\Leftrightarrow}\nc{\Impl}{\Rightarrow}
\nc\comp{^{\text{\tt c}}}\nc{\nin}{\notin}\nc{\cat}{\hat{\ }}\nc{\Un}{\bigcup}
\nc{\sub}{\subseteq}\nc{\spst}{\supseteq}\nc{\sm}{\setminus}\nc{\as}{\subseteq^*}\nc{\rest}{\restriction}
\nc{\la}{\langle}\nc{\ra}{\rangle}\nc{\dom}{\op{dom}}
\nc{\cov}{\op{cov}}\nc{\add}{\op{add}}\nc{\cof}{\op{cof}}\nc{\cf}{\op{cf}}\nc{\non}{\op{non}}\nc{\unif}{\op{non}}
\nc{\COV}{\op{COV}}\nc{\ADD}{\op{ADD}}\nc{\COF}{\op{COF}}\nc{\NON}{\op{NON}}
\newtheorem{thm}{Theorem}[section]\nc{\bthm}{\begin{thm}} \nc{\ethm}{\end{thm}}
\newtheorem{prop}[thm]{Proposition}\nc{\bprp}{\begin{prop}} \nc{\eprp}{\end{prop}}
\newtheorem{fact}[thm]{Fact}\nc{\bfct}{\begin{fact}} \nc{\efct}{\end{fact}}
\newtheorem{prob}[thm]{Problem}\nc{\bprb}{\begin{prob}} \nc{\eprb}{\end{prob}}
\newtheorem{lem}[thm]{Lemma}\nc{\blem}{\begin{lem}} \nc{\elem}{\end{lem}}
\newtheorem{claim}[thm]{Claim}\nc{\bclm}{\begin{claim}} \nc{\eclm}{\end{claim}}
\newtheorem{cor}[thm]{Corollary}\nc{\bcor}{\begin{cor}} \nc{\ecor}{\end{cor}}
\newtheorem{conj}[thm]{Conjecture}\nc{\bcnj}{\begin{conj}} \nc{\ecnj}{\end{conj}}
\theoremstyle{definition}\newtheorem{defn}[thm]{Definition}\nc{\bdfn}{\begin{defn}} \nc{\edfn}{\end{defn}}
\theoremstyle{remark}
\newtheorem{rem}[thm]{Remark}\nc{\brem}{\begin{rem}}\nc{\erem}{\end{rem}}
\newtheorem{cnv}[thm]{Convention}\nc{\bcnv}{\begin{cnv}} \nc{\ecnv}{\end{cnv}}
\newtheorem{exam}[thm]{Example}\nc{\bexm}{\begin{exam}}\nc{\eexm}{\end{exam}}
\nc{\bpf}{\begin{proof}}\nc{\epf}{\end{proof}}\nc{\be}{\begin{enumerate}}\nc{\ee}{\end{enumerate}}
\nc{\bi}{\begin{itemize}}\nc{\itm}{\item}\nc{\ei}{\end{itemize}}
\nc{\sone}{\mathsf{S}_1}\nc{\sfin}{\mathsf{S}_\mathrm{fin}}\nc{\ufin}{\mathsf{U}_\mathrm{fin}}\nc{\Sel}{\mathsf{S}}
\nc{\Split}{\mathsf{Split}}\nc{\gone}{\mathsf{G}_1}\nc{\gfin}{\mathsf{G}_\mathrm{fin}}

\title{Diagonalizations of dense families}

\author[M. Bonanzinga]{Maddalena Bonanzinga}
\author[F. Cammaroto]{Filippo Cammaroto}
\author[B. Pansera]{Bruno Antonio Pansera}
\address[Bonanzinga, Cmmaroto, Pansera]{Dipartimento di Matematica, Universit\'a di Messina, 98166 Messina, Italy}
\email{[mbonanzinga, camfil,  bpansera]@unime.it}

\author[B. Tsaban]{Boaz Tsaban}
\address[Tsaban]{Department of Mathematics, Bar-Ilan University, Ramat-Gan 52900, Israel}
\email{tsaban@math.biu.ac.il}
\urladdr{http://www.cs.biu.ac.il/\~{}tsaban}

\keywords{
dense families,
selection principles,
$\sfin(\Op,\Op)$,
Hurewicz property, 
Menger property,
$\sone(\Op,\Op)$,
C$''$, 
Rothberger property,
$\sfin(\D,\D)$,
$\sfin(\fD,\fD)$,
$\sfin(\cD,\cD)$,
M-separable,
selectively separable,
SS,
$\sone(\D,\D)$,
$\sone(\fD,\fD)$,
$\sone(\cD,\cD)$,
R-separable,
$\sfin(\Do,\D)$,
tiny sequence,
$\sfin(\cD,\cD)$,
$\sone(\Do,\D)$,
$\sone(\cD,\cD)$,
1-tiny sequence,
selectively c.c.c.,
$\sfin(\Op,\D)$,
weakly Hurewicz,
weakly Menger,
$\sone(\Op,\D)$,
weakly C$''$,
$\sone(\cO,\cD)$,
weakly Rothberger.
}

\subjclass{
Primary: 37F20; 
Secondary 26A03, 
03E75 
03E17 
}

\begin{document}

\begin{abstract}
We develop a unified framework for the study of classic and new properties involving
diagonalizations of dense families in topological spaces. We provide complete
classification of these properties. Our classification draws upon a large number of methods and
constructions scattered in the literature, and on some novel results concerning the classic properties.
\end{abstract}

\maketitle

\section{Introduction}

The following diagonalization prototypes are ubiquitous in the mathematical literature (see, e.g., the
surveys \cite{LecceSurvey, KocSurv, ict}):
\begin{description}
\item[$\sone(\scrA,\scrB)$] For all $\cU_1,\cU_2,\dots\in\scrA$, there are
$U_1\in\cU_1,U_2\in\cU_2,\dots$ such that $\sseq{U_n}\in\scrB$.
\item[$\sfin(\scrA,\scrB)$] For all $\cU_1,\cU_2,\dots\in\scrA$, there are
finite $\cF_1\sub\cU_1,\cF_2\sub\cU_2,\dots$ such that $\Un_n\cF_n\in\scrB$.
\end{description}
The papers \cite{coc1,coc2} have initiated the simultaneous consideration of these properties
in the case where $\scrA$ and $\scrB$ are important
families of open covers of a topological space $X$.
This unified study of topological properties, that were previously studied separately, had tremendous success,
some of which surveyed in the above-mentioned surveys. The field of \emph{selection principles}
is growing rapidly, and dozens of new papers appeared since these survey articles were published.
The purpose of the present paper is to initiate a similar program
for the case where $\scrA$ and $\scrB$ are dense families, as we now define.

\bdfn
Let $X$ be a topological space.
A family $\cU\sub P(X)$ is a \emph{dense family} if $\Un\cU$ is a dense subset of $X$.
A family $\cU\sub P(X)$ is in:
\begin{description}
\itm[$\D$] if $\cU$ is dense;
\itm[$\Do$] if $\cU$ is dense and all members of $\cU$ are open; and
\itm[$\Op$] if $\cU$ is an open cover of $X$.
\end{description}
\edfn
In other words, $\cU$ is a dense family if each open set in $X$ intersects some member of $\cU$.
Note that
$$\Op\sub\Do\sub\D.$$
Every element of $\D$ is refined by a dense family of singletons. It follows, for example,
that $\sfin(\D,\D)$ is equivalent to the following property, studied under various names
in the literature (see Table \ref{table} below):
\begin{quote}
For each sequence $A_n$, $n\in\N$, such that $\cl{A_n}=X$ for all $n$, there are finite
sets $F_n\sub A_n$, $n\in\N$, such that $\cl{\Un_nF_n}=X$.
\end{quote}
We study all properties $\Sel(\scrA,\scrB)$ for
$\Sel\in\{\sone,\sfin\}$ and $\scrA,\scrB\in\{\Op,\Do,\D\}$, by making use of their inter-connections.
This approach is expected to have impact beyond these properties, not only concerning properties
that imply or are implied by the above-mentioned properties (e.g., the corresponding game-theoretic
properties), but also concerning formally unrelated properties that have similar flavor.

The properties we are studying here were studied in the literature under various,
sometimes pairwise incompatible, names.
Examples are given in Table \ref{table} below, with some references. We do not give references for
$\sfin(\Op,\Op)$ and $\sone(\Op,\Op)$, because there are hundreds of them. Instead, we refer
to the above-mentioned surveys.
In this table, by \emph{obsolete} we mean that nowadays the name stands for another property.

A topological space is \emph{$\scrA$-Lindel\"of} ($\scrA\in\{\D,\Do,\Op,\dots\}$) if each member
of $\scrA$ contains a countable member of $\scrA$. If $X$ satisfies $\sfin(\scrA,\scrA)$, then
$X$ is $\scrA$-Lindel\"of. This, $\sfin(\Op,\Op)$ spaces are Lindel\"of,
$\sfin(\D,\D)$ spaces are separable, and $\sfin(\Do,\D)$ spaces are $\Do$-Lindel\"of, or equivalently,
c.c.c.\ (i.e., such that every maximal disjoint family of open sets in the space is countable).
For the latter assertion, note that every maximal disjoint family of open sets is dense, and
every maximal open refinement of an element of $\Do$ is a maximal disjoint family of open sets.
This also explains why Aurichi's notion of \emph{selectively c.c.c.} \cite{Aurichi}
is equivalent to $\sone(\Do,\D)$.\footnote{We thank Angelo Bella for bringing Aurichi's
notion of selectively c.c.c.\ to our attention, and for pointing out its equivalence to
$\sone(\Do,\D)$.}

\nc{\is}[2]{#1} 

\begin{table}[!htp]
\caption{Earlier names of the studied properties}\label{table}
\begin{center}
\vspace{0.4cm}
\begin{tabular}{rl}
\textbf{Property} & \textbf{Classic names and references}
\\ [1ex] \hline \\ [-1.5ex]
$\sfin(\Op,\Op)$ & Hurewicz (obsolete), Menger
\\ [1ex] \hline \\ [-1.5ex]
$\sone(\Op,\Op)$ & C$''$, Rothberger
\\ [1ex] \hline \\ [-1.5ex]

$\sfin(\D,\D)$ &
$\sfin(\fD,\fD)$ \cite{\is{coc6}{99}},
$\sfin(\cD,\cD)$ \cite{\is{BabinkostovaSS}{09}},
M-separable \cite{\is{BBM}{09},\is{BBMaddendum}{10},\is{RZ}{11}}
\\
&
selectively separable (SS) \cite{\is{BBMT}{08},\is{BabinkostovaSS}{09},\is{BarmanDow}{11},
\is{GS}{11},\is{BarmanDow2}{12},\is{BMS}{12},\is{Iurato}{12},\is{SakaiSS}{12},\is{BBM1}{infty},\is{BellaSSp}{infty}}
\\ [1ex] \hline \\ [-1.5ex]

$\sone(\D,\D)$ &
$\sone(\fD,\fD)$ \cite{\is{coc6}{99}},
$\sone(\cD,\cD)$ \cite{\is{BabinkostovaSS}{09}},
R-separable \cite{\is{BBM}{09},\is{BBMaddendum}{10},\is{GS}{11},\is{BMS}{12},\is{SakaiSS}{12}}
\\ [1ex] \hline \\ [-1.5ex]

$\sfin(\Do,\D)$ &
no tiny sequence
\cite{\is{Szymanski}{84},\is{Kucharski}{infty}},
$\sfin(\cD,\cD)$ \cite{\is{coc4}{98},\is{coc5}{00},\is{SakaiSS}{12}}
\\ [1ex] \hline \\ [-1.5ex]

$\sone(\Do,\D)$ &
$\sone(\cD,\cD)$ \cite{\is{coc4}{98},\is{coc5}{00},\is{SakaiSS}{12}},
no 1-tiny sequence \cite{\is{Kucharski}{infty}},
selectively c.c.c.\ \cite{\is{Aurichi}{infty}}
\\ [1ex] \hline \\ [-1.5ex]

$\sfin(\Op,\D)$ &
weakly Hurewicz (obsolete) \cite{\is{Daniels}{88},\is{SakaiSS}{12}},
weakly Menger \cite{\is{Pansera}{11}, \is{SakaiLec}{oo}}
\\ [1ex] \hline \\ [-1.5ex]

$\sone(\Op,\D)$ &
weakly C$''$ \cite{\is{Daniels}{88},\is{SakaiSS}{12}},
$\sone(\cO,\cD)$ \cite{\is{coc4}{98}},
weakly Rothberger \cite{\is{Pansera}{11}}\\
\hline
\end{tabular}
\vspace{0.2cm}
\end{center}
\end{table}

\section{Classification}

\subsection{Implications}

We begin with the 18 properties of the form $\Sel(\scrA,\scrB)$, where
$\Sel\in\{\sone,\sfin\}$ and $\scrA,\scrB\in\{\Op,\Do,\D\}$.
We first observe that six of these properties are void, and consequently need not be considered.

\blem
Let $X$ be a nondiscrete Hausdorff space. Then $X$ does not satisfy any
of the properties
$\Sel(\Do,\Op)$,
$\Sel(\D,\Op)$,
$\Sel(\D,\Do)$
($\Sel\in\{\sone,\sfin\}$).
\elem
\bpf
Note that each of these properties $\Sel(\scrA,\scrB)$ implies that each member of
$\scrA$ contains a countable member of $\scrB$.
Let $x$ be a nonisolated point.

Since $\{X\sm\{x\}\}\in\Do\sm\Op$, $X$ does not satisfy $\sfin(\Do,\Op)$.

Let $y\in X\sm\{x\}$, and let $U,V$ be disjoint open neighborhoods
of $x,y$ respectively. Then the set $U\sm\{x\}$ is not closed and not dense.
Then $\cU:=\{U\sm\{x\},(U\sm\{x\})\comp\}\in\D$, and each family of open sets
contained in $\cU$ contains at most $U\sm\{x\}$, which is not dense.
Thus, $X$ does not satisfy $\sfin(\D,\Do)$.

The other assertions follow.
\epf

The following immediate equivalences ($\Sel\in\{\sone,\sfin\}$)
eliminate the need to consider 4 additional properties:
\begin{eqnarray*}
\Sel(\Do,\D) & = & \Sel(\Do,\Do)\\
\Sel(\Op,\D) & = & \Sel(\Op,\Do)
\end{eqnarray*}
We are thus left with the following eight properties.
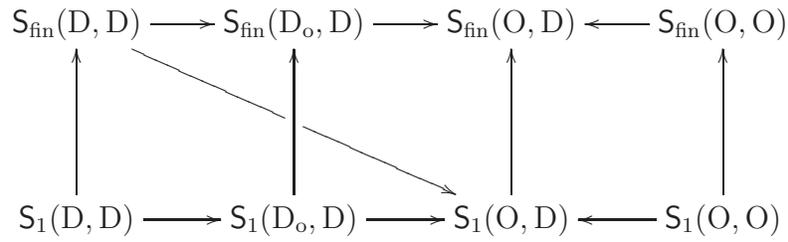
\begin{figure}[!htp]
$$\xymatrix{
\sfin(\D,\D) \ar[r]\ar'[rd][rrdd]& \sfin(\Do,\D) \ar[r]& \sfin(\Op,\D) & \sfin(\Op,\Op)\ar[l]\\
& & \\
\sone(\D,\D) \ar[r]\ar[uu]& \sone(\Do,\D) \ar[r]\ar[uu]& \sone(\Op,\D)\ar[uu] & \sone(\Op,\Op)\ar[l]\ar[uu]
}$$
\caption{The Dense Families Diagram}\label{DFD}
\end{figure}
For the diagonal implication, note that $\sfin(\D,\D)$ implies separability, which in turn implies $\sone(\Op,\D)$.
Indeed, every countable space satisfies $\sone(\Op,\Op)$, and we have
the following.

\bprp
Let $\Sel\in\{\sone,\sfin\}$. If $X$ has a dense subset satisfying $\Sel(\Op,\Op)$,
then $X$ satisfies $\Sel(\Op,\D)$. \qed
\eprp

\forget
\brem
Angelo Bella pointed out to us that Aurichi introduces in \cite{Aurichi}
a property like $\sone(\Do,\D)$, with the formal difference that maximal antichains
of open sets are given, instead of dense families of open sets. However, Bella points out that,
since a maximal antichain refinement of an element of $\Do$ is also in $\Do$,
Aurichi's property is equivalent to $\sone(\Do,\D)$.
\erem
\forgotten

\subsection{Non-implications}

To make it clear which properties are possessed by the examples given below and which not,
we supply a version of the Dense Families Diagram (Figure \ref{DFD}) with a full bullet ($\bullet$)
for each property the example satisfies, and an empty bullet ($\circ$) for each property not satisfied by the example.

\nc{\Setting}[8]{
\xymatrix{
#1 \ar[r]\ar'[rd][rrdd]& #2 \ar[r]& #3 & #4\ar[l]\\
& & \\
#5 \ar[r]\ar[uu]& #6 \ar[r]\ar[uu]& #7 \ar[uu] & #8\ar[l]\ar[uu]
}}

\subsubsection{Uncountable examples}

\bprp\label{betaN}\label{Baire}
\mbox{}
\be
\item
The spaces $\R$ and $\beta\N$ satisfy the following setting.
$$\Setting{\bullet}{\bullet}{\bullet}{\bullet}{\bullet}{\bullet}{\bullet}{\circ}$$

\item
The Baire space $\NN$ satisfies the following setting.
$$\Setting{\bullet}{\bullet}{\bullet}{\circ}{\bullet}{\bullet}{\bullet}{\circ}$$
\ee
\eprp
\bpf
Each of these spaces has a countable pseudobase, and thus \cite{BBMT} satisfies $\sone(\D,\D)$.

Being $\sigma$-compact, $\R$ and $\beta\N$ satisfy $\sfin(\Op,\Op)$. Since $\sone(\Op,\Op)$ subsets of $\R$ have
measure zero, $\R$ does not satisfy $\sone(\Op,\Op)$. For the same reason, the unit interval
$[0,1]$ does not satisfy $\sone(\Op,\Op)$. Since $[0,1]$ (being separable and compact) is a continuous
image of $\beta\N$ and $\sone(\Op,\Op)$ is preserved by continuous images, $\beta\N$ does not satisfy
$\sone(\Op,\Op)$, too.

The Baire space does not satisfy $\sfin(\Op,\Op)$ (e.g., \cite{coc2}).
\epf

\bthm\label{betaN-N}
$\beta\N\sm\N$ satisfies the following setting.
$$\Setting{\circ}{\circ}{\bullet}{\bullet}{\circ}{\circ}{\circ}{\circ}$$
\ethm
\bpf
Being compact, $\beta\N\sm\N$ satisfies $\sfin(\Op,\Op)$.
$\sone(\Op,\Op)$ is preserved by countable unions, and every countable space
satisfies $\sone(\Op,\Op)$. Thus, had $\beta\N\sm\N$ satisfied $\sone(\Op,\Op)$,
so would $\beta\N$, its union with the countable set $\N$, in contradiction
to Proposition \ref{betaN}.

That $\beta\N\sm\N$ does not satisfy $\sfin(\Do,\D)$ follows from the following.
\blem
$\beta\N\sm\N$ is not $\Do$-Lindel\"of, that
is, there is an element of $\Do$ with no countable subset in $\Do$.
\elem
\bpf
Let $\roth$ be the family of all infinite subsets of $\N$.
For $A\in\roth$, let
$$[A]=\{p\in\beta\N\sm\N : A\in p\}$$
be the standard basic clopen subset of $\beta\N\sm\N$.

Let $\cA\sub\roth$ be a maximal almost disjoint family.
Take $\cU=\{[A] : A\in\cA\}$.

$\cU\in\Do$: Let $[B]$ be a  basic clopen set in $\beta\N\sm\N$.
There is $A\in\cA$ such that $B\cap A$ is infinite:
Indeed, otherwise $\cA\cup\{B\}$ would be an almost disjoint family extending $\cA$,
in contradiction to the maximality of $\cA$.
Thus, the set $[B]\cap[A]=[B\cap A]$ is nonempty.

$\cU$ does not contain a countable element from $\Do$: Let $\cV=\{[A_n] : n\in\N\}$
be a countable subset of $\cU$. Since $\cA$ is a maximal almost disjoint family, $\cA$ is uncountable.
Let $B\in\cA\sm\sseq{A_n}$.
For each $n$, since $B,A_n\in\cA$ and are distinct, $B\cap A_n$ is finite.
Thus,
$$[B]\cap[A_n]=[B\cap A_n]=\emptyset.\qedhere$$
\epf
To finish the proof, we prove the following.
\blem
$\beta\N\sm\N$ does not satisfy $\sone(\Op,\D)$.
\elem
\bpf
By induction on $n$, choose for each sequence $(s_1,\dots,s_n)\in\{0,1\}^n$
an infinite set $I_{s_1,\dots,s_n}$ such that:
\be
\itm $I_0\cup I_1=\N$, and the union is disjoint.
\itm For each $n$ and each sequence $(s_1,\dots,s_n)\in\{0,1\}^n$,
$$I_{s_1,\dots,s_n,0}\cup I_{s_1,\dots,s_n,1}=I_{s_1,\dots,s_n},$$
and this union is disjoint.
\ee
For each $n$, let
$$\cU_n=\{[I_{s_1,\dots,s_n}] : s_1,\dots,s_n\in\{0,1\}\}.$$
$\cU_n\in\Op$. Indeed, since the involved unions are finite,
$$\Un_{(s_1,\dots,s_n)\in\{0,1\}^n}[I_{s_1,\dots,s_n}]=
\left[\Un_{(s_1,\dots,s_n)\in\{0,1\}^n}I_{s_1,\dots,s_n}\right]=[\N]=\beta\N\sm\N.$$
Now, consider any selection $[I_{s^n_1,\dots,s^n_n}]\in\cU_n$, $n\in\N$.
By induction on $n$, choose $t_n\in\{0,1\}$ such that
$$I_{t_1,\dots,t_n}\cap \left( I_{s^1_1}\cup I_{s^2_1,s^2_2}\cup \cdots \cup I_{s^n_1,\dots,s^n_n}\right)=\emptyset$$
for all $n$. This is possible, since the latter union is contained in a union of at most
$2^{n-1}+2^{n-2}+\cdots+2+1=2^n-1$
sets of the form $I_{s_1,\dots,s_n}$.
The sets
$I_{t_1,\dots,t_n}$, $n\in\N$, form a decreasing sequence of infinite subsets of $\N$.
Let $A$ be a pseudointersection of these sets.
For each $n$,
$$A\as I_{t_1,\dots,t_n}\sub \left( I_{s^1_1}\cup I_{s^2_1,s^2_2}\cup \cdots \cup I_{s^n_1,\dots,s^n_n}\right)\comp
\sub I_{s^n_1,\dots,s^n_n}\comp,$$
and thus $A\cap I_{s^n_1,\dots,s^n_n}$ is finite. Therefore,
$[A]\cap [I_{s^n_1,\dots,s^n_n}]=[A\cap I_{s^n_1,\dots,s^n_n}]=\emptyset$
for all $n$.
\epf
This completes the proof of Theorem \ref{betaN-N}.
\epf

The properties $\sone(\Op,\Op)$ and $\sfin(\Op,\Op)$ are hereditary for closed subsets (e.g., \cite{coc2}).
$\sfin(\Op,\D)$ is hereditary for compact subsets since compact spaces satisfy $\sfin(\Op,\Op)$.
In contrast to that, we have the following.

\bcor
None of the properties
$\sone(\D,\D)$, $\sone(\Do,\D)$, $\sfin(\D,\D)$, $\sfin(\Do,\D)$, $\sone(\Op,\D)$,
is hereditary for compact subsets.
\ecor
\bpf
Proposition \ref{betaN} and Theorem \ref{betaN-N}.
\epf

We consider ordinals $\alpha$ with the order topology, so that the basic clopen sets are the intervals
$(\beta,\gamma)$ or $[0,\beta)$ or $(\beta,\alpha)$, where $\beta,\gamma\in\alpha$.

\bthm\label{successor}
Each uncountable successor ordinal $\alpha+1$ satisfies the following setting.
$$\Setting{\circ}{\circ}{\bullet}{\bullet}{\circ}{\circ}{\bullet}{\bullet}$$
\ethm
\bpf
The theorem follows from the following two lemmata.
\blem[\cite{BCKM}]
Each uncountable successor ordinal $\alpha+1$ satisfies $\sone(\Op,\Op)$.
\elem
\bpf
For completeness, we reproduce the proof: Given open covers $\cU_1,\cU_2,\dots$ of
$\alpha+1$, each consisting of basic clopen sets,
pick $U_1\in\cU_1$ with $\alpha\in U_1$. If $\alpha\sm U_1\neq \emptyset$,
it is a successor ordinal, and we can cover its last element by some $U_2\in\cU_2$.
This must end after finitely many steps, since the sequence of last elements is decreasing.
\epf

\blem
For each uncountable ordinal $\alpha$, $\alpha+1$ is not $\Do$-Lindel\"of.
\elem
\bpf
Let $S$ be the set of all successor ordinals in $\alpha+1$.
$S$ is dense in $\alpha+1$.
Since successor ordinals are isolated, $\cU=\{\{\beta\} : \beta\in S\}\in\Do$.
Let $\cV$ be a countable subset of $\cU$. Then $\beta=\sup\{\beta\in\omega_1\cap\Un\cV\}<\omega_1$.
Thus, no element of $\cV$ intersects $(\beta,\omega_1)$.
\epf
This completes the proof of Theorem \ref{successor}.
\epf

Consider the following construction from \cite{coc4}.
The \emph{Alexandroff double} of $[0,1]$ is the space $[0,1]\x\{0,1\}$, with the basic open sets
$\{(x,1)\}$ for each $x\in [0,1]$, and $(U\x\{0,1\})\sm (F\x\{1\})$ for each open $U$ in $[0,1]$ and each
finite $F\sub [0,1]$. For each dense $X\sub [0,1]$, the subspace
$$T(X):=([0,1]\x\{0\})\cup (X\x\{1\})$$
is compact Hausdorff.\footnote{The notation $T(X)$ for this construction is due to Scheepers,
in recognition of the inspiration provided by a related construction of Tkachuk.}

\bthm\label{TX}
For a dense $X\sub[0,1]$, $T(X)$ satisfies the following:
\be
\item The setting in Proposition \ref{betaN}(1) if $X$ is countable;
\item The setting in Theorem \ref{betaN-N} if $X$ is uncountable and does not have strong measure zero; and
\item The following setting if $X$ is uncountable and has strong measure zero (e.g., when $X$ is
a Luzin set).
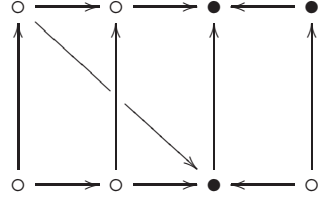
\ee
$$\Setting\circ\circ\bullet\bullet\circ\circ\bullet\circ$$
\ethm
\bpf
The unit interval $[0,1]$ is a closed subspace of $T(X)$.
Since $\sone(\Op,\Op)$ is preserved by moving to closed subsets, $T(X)$ does not satisfy $\sone(\Op,\Op)$.
On the other hand, $T(X)$ is compact (after covering its lower part by finitely many sets, there remain only
finitely many uncovered points on its top part), and thus satisfies $\sfin(\Op,\Op)$.

(1) In this case, $T(X)$ has a countable base, and thus satisfies $\sone(\D,\D)$.

(2,3) If $X$ is uncountable, then $T(X)$ is not $\Do$-Lindel\"of.
Indeed, $X\x\{1\}$ is an uncountable discrete dense subspace of $T(X)$.
Since being $\Do$-Lindel\"of is hereditary for dense subspaces, had $T(X)$ been $\Do$-Lindel\"of,
so would the uncountable discrete space $X\x\{1\}$, a contradiction.
It remains to consider $\sone(\Op,\D)$, and this was done by Scheepers, who proved in \cite{coc4}
that $T(X)$ has this property if and only if $X$ has strong measure zero.
\epf

Daniels \cite{Daniels} proved that, for each $\Sel\in\{\sone,\sfin\}$,
if every finite subproduct of a product space $\prod_{i\in I}X_i$ satisfies
$\Sel(\Op,\D)$, then so does the full product $\prod_{i\in I}X_i$.
It is a classic fact that the same assertion holds for \emph{$\Do$-Lindel\"of} (equivalently, c.c.c.)
spaces.
We prove that this is also the case for $\Sel(\Do,\D)$ ($\Sel\in\{\sone,\sfin\}$).
Note that this is \emph{not} the case for the remaining properties: Consider
the countably infinite power $\NN$ of $\N$ for $\Sel(\Op,\Op)$
and the (nonseparable) power $\N^{\aleph_1}$ for $\Sel(\D,\D)$.

Modulo Lemma \ref{finprod} below, whose proof is similar to that of Theorem 2.27 in \cite{BBMT},
the following Theorem \ref{prodthm} was independently proved by Leandro Aurichi \cite{Aurichi}.

\bthm\label{prodthm}
Let $\Sel\in\{\sone,\sfin\}$.
Let $X_i$, $i\in I$, be spaces such that $\prod_{i\in F}X_i$ satisfies
$\Sel(\Do,\D)$ for all finite $F\sub I$. Then $\prod_{i\in I}X_i$ satisfies
$\Sel(\Do,\D)$.
\ethm
\bpf
We prove the assertion for $\Sel=\sone$. The proof in the other case is similar.

For an open set $U$ in a product space $\prod_iX_i$, let $\supp(U)$ (the \emph{support} of $U$)
be the finite set of coordinates $i$ where $\pi_i(U)\neq X_i$.
Note that open sets $U,V$ in a product space intersect if and only if their
projections $\pi_F(U),\pi_F(V)$ intersect, for $F=\supp(U)\cap\supp(V)$.
\blem\label{finprod}
Let $X_n$, $n\in \N$, be spaces such that $\prod_{n\le k} X_n$ satisfies
$\sone(\Do,\D)$ for all $k$. Then $\prod_{n\in \N}X_n$ satisfies
$\sone(\Do,\D)$.
\elem
\bpf
Let $X=\prod_{n\in\N}X_n$, and let $\cU_1,\cU_2,\dots\in\Do(X)$.
Decompose $\N=\Un_{k\in\N}A_k$, with each $A_k$ infinite.

Fix $k\in\N$.
Since $\prod_{n\le k}X_n$ satisfies $\sone(\Do,\D)$ and
$\set{\pi_{\{1,\dots,k\}}(U)}{U\in\cU_n}$ is in $\Do(\prod_{n\le k}X_n)$
for all $n\in A_k$, there are $U_n\in\cU_n$, $n\in A_k$, such that
$\set{\pi_{\{1,\dots,k\}}(U_n)}{n\in A_k}\in\D(\prod_{n\le k}X_n)$.

We claim that $\set{U_n}{n\in\N}$ is a dense family in $X$. Indeed, let $U$ be an open
subset of $X$. Let $k$ be such that $\supp(U)\sub\{1,\dots,k\}$.
By our construction, there is $n\in A_k$ such that the projections $\pi_{\{1,\dots,k\}}(U)$ and $\pi_{\{1,\dots,k\}}(U_n)$
intersect. Since $\supp(U)\sub\{1,\dots,k\}$, $U$ intersects $U_n$.
\epf

We now prove the general assertion.
Let $X=\prod_{i\in I}X_i$, and  let $\cU_1,\cU_2,\dots\in\Do(X)$.
Decompose $\N=\Un_{k\in\N}A_k$, with each $A_k$ infinite.

Let $I_1$ be any countable nonempty subset of $I$.
By the lemma, $\prod_{i\in I_1}X_i$ satisfies $\sone(\Do,\D)$.
Thus, there are $U_n\in\cU_n$, $n\in A_1$, such that
$\set{\pi_{I_1}(U_n)}{n\in A_1}\in\D(\prod_{i\in I_1}X_i)$.
Let
$$I_2=I_1\cup\Un_{n\in A_1}\supp(U_n),$$
and note that $I_2$ is countable.
By the lemma, $\prod_{i\in I_2}X_i$ satisfies $\sone(\Do,\D)$.
Thus, there are $U_n\in\cU_n$, $n\in A_2$, such that
$\set{\pi_{I_2}(U_n)}{n\in A_2}\in\D(\prod_{i\in I_2}X_i)$.
Let
$$I_3=I_2\cup\Un_{n\in A_2}\supp(U_n).$$
Continue in the same manner.

We claim that $\set{U_n}{n\in\N}$ is a dense family in $X$. Indeed, let $U$ be an open
subset of $X$. Let $I_\infty=\Un_{n\in\N}I_n$, and $F=\supp(U)\cap I_\infty$.
Let $k$ be such that $F\sub I_k$.
By the construction, $\pi_F(U)$ intersects some $\pi_F(U_n)$, $n\in A_k$.
Since $\supp(U_n)\sub I_\infty$, $\supp(U)\cap\supp(U_n)\sub F$.
Thus, $U$ intersects $U_n$.
\epf

\bthm\label{Rpow}
For each nonempty set $X$, the Tychonoff power $\R^X$ satisfies:
\be
\itm The setting of Proposition \ref{betaN}(1) if $X$ is finite nonempty;
\itm The setting of Proposition \ref{Baire}(2) if $X$ is countably infinite; and
\itm The following setting if $X$ is uncountable.
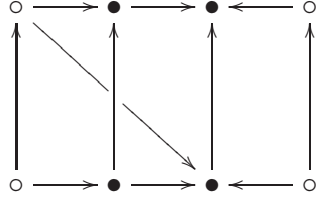
\ee
$$\Setting\circ\bullet\bullet\circ\circ\bullet\bullet\circ$$
\ethm
\bpf
If $X$ is countable, then $\R^X$ has a countable base, and thus satisfies $\sone(\D,\D)$.
If $X$ is finite, then $\R^X$ is $\sigma$-compact, and thus satisfies $\sfin(\Op,\Op)$.
Since $\R$ is a continuous image of $\R^X$, $\R^X$ does not satisfy $\sone(\Op,\Op)$.
This concludes (1).

(2) $\R^\N$ does not satisfy $\sfin(\Op,\Op)$.

(3) As $X$ is uncountable, $\R^X$ is neither separable nor Lindel\"of,
and in particular does not satisfy $\sfin(\D,\D)$ or $\sfin(\Op,\Op)$.
It remains to prove
that $\R^X$ satisfies $\sone(\Do,\D)$, and this follows from Theorem \ref{prodthm}.\footnote{That
$\R^X$ satisfies $\sone(\Do,\D)$ was also, independently, proved by Aurichi \cite{Aurichi}.}
\epf

\bthm\label{2pow}
The Tychonoff power $\{0,1\}^X$ satisfies:
\be
\itm The setting of Proposition \ref{betaN}(1) if $X$ is countably infinite; and
\itm The following setting if $X$ is uncountable.
$$\Setting\circ\bullet\bullet\bullet\circ\bullet\bullet\circ$$
\ee
\ethm
\bpf
$\sfin(\Op,\Op)$ for $\{0,1\}^X$ follows from compactness.

The Cantor space does not satisfy $\sone(\Op,\Op)$,
e.g., since $[0,1]$ is its continuous image.
Thus, $\{0,1\}^X$ does not satisfy $\sone(\Op,\Op)$.

By Theorem \ref{prodthm}, $\{0,1\}^X$ satisfies $\sone(\Do,\D)$.

Finally, if $X$ is uncountable, then $\{0,1\}^X$ is not separable, and thus does not satisfy $\sfin(\D,\D)$.
\epf

For a topological space $X$, let $\Om$ be the family of all $\cU\in\Op$ such
that every finite subset of $X$ is contained in some member of $\cU$, and $X\notin\cU$.
Covering properties involving this family were studied extensively \cite{LecceSurvey, KocSurv, ict}.

\bthm\label{CpX}
Let $X$ be an infinite Tychonoff space. The space $C_p(X)$ satisfies:
\be
\itm The setting
$$\Setting\circ\bullet\bullet\circ\circ\bullet\bullet\circ$$
if $X$ does not satisfy $\sfin(\Omega,\Omega)$ (e.g., if $X=\NN$) or there is no coarser, second countable Tychonoff
topology on $X$;
\itm The setting
$$\Setting\bullet\bullet\bullet\circ\circ\bullet\bullet\circ$$
if $X$ satisfies $\sfin(\Omega,\Omega)$ but not $\sone(\Omega,\Omega)$
and there is a coarser, second countable Tychonoff
topology on $X$ (e.g., if $X=\R$); and
\itm The setting of Proposition \ref{Baire}(2)
if $X$ satisfies $\sone(\Omega,\Omega)$
and there is a coarser, second countable Tychonoff topology on $X$.
\ee
\ethm
\bpf
$C_p(X)$ is dense in $\R^X$. By Theorem \ref{Rpow}, $C_p(X)$ satisfies $\sone(\Do,\D)$.
As $X$ is infinite, $C_p(X)$ does not satisfy $\sfin(\Op,\Op)$ \cite[Theorem II.2.10]{Arhan92}.

By \cite[Theorems 21,57]{BBM}, $C_p(X)$ satisfies $\sfin(\D,\D)$ (respectively, $\sone(\D,\D)$) if and only if
there is a coarser, second countable Tychonoff topology on $X$ and $X$ satisfies $\sfin(\Om,\Om)$
(respectively, $\sone(\Om,\Om)$).
\epf

Let $\R_\mathrm{coc}$ be $\R$ with the topology generated by the usual open intervals and all cocountable
sets. This example was first considered in this context by Aurichi \cite{Aurichi}.

\bprp\label{Rcoc}
The space $\R_\mathrm{coc}$ satisfies the setting of Theorem \ref{2pow}(2).
\eprp
\bpf
$\R_\mathrm{coc}$ is not separable, and thus does
not satisfy $\sfin(\D,\D)$.
Aurichi \cite{Aurichi} proved that $\R_\mathrm{coc}$ satisfies $\sone(\Do,\D)$.

$\R_\mathrm{coc}$ does not satisfy $\sone(\Op,\Op)$ because $\R$, which is coarser,
does not.

$\R_\mathrm{coc}$ satisfies $\sfin(\Op,\Op)$ because $\R$ does: Given $\cU_1,\cU_2,\dots\in\Op(\R_\mathrm{coc})$,
whose elements have the form $(a,b)\sm C$ with $C$ countable,
let $\cU_1',\cU_2',\dots$ be the open covers of $\R$ obtained by replacing each $(a,b)\sm C$ with $(a,b)$.
Take finite $\cF_1'\sub\cU_1,\cF_2'\sub\cU_2',\dots$ such that $\Un_n\cF_n'$ is a cover of $\R$.
Then moving back to the original elements, we have that $\R\sm\Un_n\cF_n$ is countable. Choose one more
element from each $\cU_n$ to cover this countable remainder.
\epf

Let $X$ be a topological space. The \emph{Pixley--Roy space} $\PR(X)$ is the space of all nonempty finite
subsets of $X$, with the topology determined by the basic open sets
$$[F,U] := \set{H\in\PR(X)}{F\sub H\sub U},$$
$F\in\PR(X)$ and $U$ open in $X$.

For regular spaces $X$, the Pixley--Roy space $\PR(X)$ is zero-dimensional, completely regular,
and hereditarily metacompact.

\bthm\label{PRunc}
Let $X$ be an uncountable separable metrizable space. The Pixley--Roy space $\PR(X)$ satisfies:
\be
\itm The setting of Theorem \ref{CpX}(1) if $X$ satisfies $\sone(\Om,\Om)$;
\itm None of the properties if $X$ does not satisfy $\sfin(\Om,\Om)$; and
\itm The following setting if $X$ satisfies $\sfin(\Om,\Om)$ but not $\sone(\Om,\Om)$
(e.g., if $X=\R$).
\ee
$$\Setting\circ\bullet\bullet\circ\circ\circ\circ\circ$$
\ethm
\bpf
Daniels \cite{Daniels} proved that, for a metrizable space $X$, $\PR(X)$ satisfies $\sone(\Op,\D)$ (respectively,
$\sfin(\Op,\D)$) if and only if $X$ satisfies $\sone(\Om,\Om)$ (respectively, $\sfin(\Om,\Om)$).

Scheepers proved that, for Pixley--Roy spaces of separable metrizable spaces
and $\Sel\in\{\sone,\sfin\}$,
$\Sel(\Do,\D)=\Sel(\Op,\D)$ \cite{coc5}.

If $\PR(X)$ satisfies $\sfin(\D,\D)$, then it is separable. It is a classic fact that, in this case, $X$ is countable
(references are available in \cite{SakaiSS}). Thus, in our case, $\PR(X)$ does not satisfy $\sfin(\D,\D)$.

\blem
The following are equivalent, for a topological space $X$:
\be
\itm $\PR(X)$ has a countable cover by basic open sets;
\itm $\PR(X)$ is Lindel\"of;
\itm $\PR(X)$ satisfies $\sfin(\Op,\Op)$;
\itm $\PR(X)$ satisfies $\sone(\Op,\Op)$;
\itm $\PR(X)$ is countable; and
\itm $X$ is countable.
\ee
\elem
\bpf[Proof of $(1)\Impl(6)$]
Assume that $\PR(X)=\Un_n[F_n,U_n]$.
For each $x\in\PR(X)$, let $n$ be such that $\{x\}\in [F_n,U_n]$. Then $F_n\sub\{x\}$,
that is, $F_n=\{x\}$. It follows that there are only countably many singletons in $\PR(X)$,
that is, $X$ is countable.
\epf
This completes the proof of Theorem \ref{PRunc}.
\epf

\brem
By recent results of Sakai \cite{SakaiLec}, Theorem \ref{PRunc} generalizes from separable
metrizable spaces to semi-stratifiable ones.
\erem

\subsubsection{Countable examples}

As pointed out already, Scheepers proved that, for Pixley--Roy spaces of separable metrizable spaces,
$\Sel(\Do,\D)=\Sel(\Op,\D)$ for both $\Sel\in\{\sone,\sfin\}$ \cite{coc5}.
We prove an analogous assertion for countable spaces (note the difference in the
properties involved).

\bthm\label{PRctblchar}
Let $\Sel\in\{\sone,\sfin\}$.
Let $X$ be a countable topological space. Then
$\PR(X)$ satisfies $\Sel(\D,\D)$ if and only if it satisfies $\Sel(\Do,\D)$.
\ethm
\bpf
We prove the assertion for $\Sel=\sone$. The proof of the remaining assertion is similar.

Assume that $\PR(X)$ satisfies $\sone(\Do,\D)$, and let
$D_1,D_2,\dots$ be dense subsets of $\PR(X)$.
Fix an enumeration $\PR(X)=\set{H_n}{n\in\N}$, and a partition $\N=\Un_kI_k$ with each
$I_k$ infinite.

Fix $k$. For each $n\in I_k$, the family
$$\cU_n=\set{[F,X]}{H_k\sub F\in D_n}$$
is dense open in the subspaces $[H_k,X]$ of $\PR(X)$:
For each basic open $[H,U]$ in $\PR(X)$ with
$[H\cup H_k,U]=[H,U]\cap [H_k,X]\neq\emptyset$,
let $F\in D_n\cap [H\cup H_k,U]$. Then $[F,X]\in\cU_n$,
and $F\cup H\cup H_k\in [F,X]\cap [H\cup H_k,U]$.

Since $\sone(\Do,\D)$ is hereditary for open subsets, there are for each $k$ elements
$[F_n,X]\in\cU_n$, $n\in I_k$, such that $\set{[F_n,X]}{n\in I_k}\in\D([H_k,X])$.
It remains to observe that $\set{F_n}{n\in\N}$ is dense in $\PR(X)$.
Indeed, let $[F,U]$ be a nonempty basic open set in $\PR(X)$. Let $k$ be such that
$H_k=F$. Since $[H_k,U]$ is open in $[H_k,X]$, there is $n\in I_k$ (so that $H_k\sub F_n$)
such that
$$[F_n,U] = [F_n\cup H_k,U] = [F_n,X]\cap [H_k,U]\neq \emptyset.$$
Then $F_n\in [H_k,U]$.
\epf

\nc{\PNfin}{\pi\mathrm{N}_\mathrm{fin}}
For a topological space $X$ and a point $x\in X$, $\PNfin(x)$ be the family of
all $\pi$-networks $\cN$ of $x$ (i.e., such that each open neighborhood of $x$ contains an element of $\cN$)
such that all members of $\cN$ are finite.
For $\Sel\in\{\sone,\sfin\}$, say that $X$ satisfies $\Sel(\PNfin, \PNfin)$ if
$\Sel(\PNfin(x), \PNfin(x))$ holds for all $x\in X$.

\bcor\label{PRcor}
Let $X$ be a countable topological space, and $\Sel\in\{\sone,\sfin\}$.
The following assertions are equivalent:
\be
\itm $\PR(X)$ satisfies $\Sel(\Do,\D)$;
\itm $\PR(X)$ satisfies $\Sel(\D,\D)$;
\itm $X$ is countable, and all finite powers of $X$ satisfy $\Sel(\PNfin, \PNfin)$.
\ee
\ecor
\bpf
The equivalence of (2) and (3) was proved by Sakai \cite{SakaiSS}.
Apply Theorem \ref{PRctblchar}.
\epf

\bthm\label{h}
Let $X$ be a countable topological space. The Pixley--Roy space $\PR(X)$ satisfies:
\be
\itm All properties in the diagram if all finite powers of $X$ satisfy $\sone(\PNfin,\PNfin)$;
\itm The setting
$$\Setting\bullet\bullet\bullet\bullet\circ\circ\bullet\bullet$$
if some finite power of $X$ does not satisfy $\sone(\PNfin,\PNfin)$, but all finite
powers of $X$ satisfy $\sfin(\PNfin,\PNfin)$;
\itm The setting
$$\Setting\circ\circ\bullet\bullet\circ\circ\bullet\bullet$$
if some finite power of $X$ does not satisfy $\sfin(\PNfin,\PNfin)$.
\ee
\ethm
\bpf
Since $X$ is countable, so is $\PR(X)$. Thus, $\PR(X)$ satisfies $\sone(\Op,\Op)$.
Apply Corollary \ref{PRcor}.
\epf

To obtain concrete examples from Theorem \ref{h}, we use Nyikos' \emph{Cantor Tree} topologies
and a result of Sakai.
Let $\{0,1\}^{<\oo}$ be the set of all finite sequences in $\{0,1\}$.
For $s,t\in \{0,1\}^{<\oo}$, let $s\sub t$ mean that $t$ is an end-extension of $s$.

\nc{\CantorTree}{\op{CT}}

Let $X\sub\Cantor$, and define a topology on the countable space
$\CantorTree(X):=\{0,1\}^{<\oo}\cup\{\oo\}$ by
declaring all points of $\{0,1\}^{<\oo}$ isolated, and taking the sets
$$\CantorTree(X)\sm\left(\{0,1\}^{\le k}\cup\set{s\in \{0,1\}^{<\oo}}{\exists f\in F,\ s\sub f}\right),$$
$k\in\N$, $F\sub X$ finite, as a local base at $\oo$.

\bthm[Sakai \cite{SakaiSS}]
Let $X\sub\Cantor$ and $\Sel\in\{\sone,\sfin\}$. The following assertions are equivalent:
\be
\itm $\CantorTree(X)$ satisfies $\Sel(\PNfin,\PNfin)$;
\itm $X$ satisfies $\Sel(\Om,\Om)$.
\ee
\ethm

The first construction of a countable space not satisfying $\sfin(\Do,\D)$ is
due to Aurichi \cite{Aurichi}. Our method makes it possible to transport examples
from classic selection principles, and is consequently more flexible, as the following
theorem shows.

\bthm\label{t}
Let $X\sub\Cantor$. The countable space $\PR(\CantorTree(X))$ satisfies:
\be
\itm Setting (2) in Theorem \ref{h} if $X$ satisfies $\sfin(\Om,\Om)$ but not $\sone(\Om,\Om)$ (e.g., $X=\Cantor$);
\itm Setting (3) in Theorem \ref{h} if $X$ does not satisfy $\sfin(\Om,\Om)$ (e.g., $X=\NN$). \qed
\ee
\ethm

We conclude with an example of Barman and Dow \cite{BarmanDow}: Let
$\compactN$ be the one-point compactification of $\N$. Take the box-product on $\NcompactN$.
Let
\nc{\EI}{\mathbb{EI}^\Box}
$$\EI=\set{f\in\NcompactN}{\exists n,\ f(1),\dots,f(n)<\oo,\ f(n+1)=f(n+2)=\cdots=\oo},$$
a countable subspace of the box-product space $\NcompactN$.
For each $n$, the set $A_n :=\set{f\in X}{f(1),\dots,f(n)<\oo}$
is dense in $\EI$, and for each selection of finite sets $F_n\sub A_n$,
$\EI\cap\prod_n[\max\set{f(n)}{f\in F_n}+1,\oo]$ is nonempty and disjoint from all $F_n$.
Thus, $\EI$ does not satisfy $\sfin(\D,\D)$ \cite{BarmanDow}.

\bthm
The Barman--Dow space $\EI$ satisfies the following setting.
$$\Setting\circ\bullet\bullet\bullet\circ\bullet\bullet\bullet$$
\ethm
\bpf
Since $\EI$ is countable, it satisfies $\sone(\Op,\Op)$.
It Remains to prove that $\EI$ satisfies $\sone(\Do,\D)$.
The proof is similar to that of Theorem \ref{prodthm}.

Let $\cU_1,\cU_2,\dots\in\Do(\EI)$.
Decompose $\N=\Un_{k\in\N}I_k$, with each $I_k$ infinite.

Fix $k\in\N$.
Since $\pi_{\{1,\dots,k\}}(\EI)$ has a countable base,
it satisfies $\sone(\Do,\D)$. As the family
$\set{\pi_{\{1,\dots,k\}}(U)}{U\in\cU_n}$ is dense open in $\pi_{\{1,\dots,k\}}(\EI)$
for all $n\in I_k$, there are $U_n\in\cU_n$, $n\in I_k$, such that
$\set{\pi_{\{1,\dots,k\}}(U_n)}{n\in I_k}$ is dense in $\pi_{\{1,\dots,k\}}(\EI)$.

We claim that $\set{U_n}{n\in\N}$ is a dense family in $\EI$.
For an open set in $\NcompactN$ intersecting $\EI$, let $F(U)=\set{k}{\oo\notin\pi_k(U)}$.
Then $F(U)$ is finite.

Let $k$ be such that $F(U)\sub\{1,\dots,k\}$.
Let $n\in I_k$ be such that the projections $\pi_{\{1,\dots,k\}}(U)$ and $\pi_{\{1,\dots,k\}}(U_n)$
intersect. Since $F(U)\sub\{1,\dots,k\}$, $U$ intersects $U_n$.
\epf

\bthm
No implication can be added to the
Dense Families Diagram (Figure \ref{DFD}), except for those obtained by composition of existing ones.
Moreover, this is exhibited by ZFC examples.
\ethm
\bpf
We go over the properties one by one, and verify that no new implication can be added from it
to another property, by referring to the an appropriate (one, in case there are several)
number of a proposition or a theorem.
When treating a property, we consider only potential implications not ruled
out by the treatment of the previous properties.
\be
\itm $\sone(\D,\D)\nrightarrow\sfin(\Op,\Op)$ (\ref{betaN}).
\itm $\sfin(\D,\D)\nrightarrow\sfin(\Do,\D)$ (\ref{t}).
\itm $\sone(\Do,\D)\nrightarrow\sfin(\D,\D)$ (\ref{Rpow}).
\itm $\sfin(\Do,\D)\nrightarrow\sone(\Op,\D)$ (\ref{PRunc}).
\itm $\sone(\Op,\Op)\nrightarrow\sfin(\Do,\D)$ (\ref{successor}).
\itm $\sfin(\Op,\Op)\nrightarrow\sone(\Op,\D)$ (\ref{betaN-N}).\qedhere
\ee
\epf

The classification is completed.

\ed